\documentclass[]{article}
\usepackage{amssymb}
\usepackage{amsmath}
\usepackage{amssymb}
\usepackage{tikz}
\usepackage{graphicx,epsfig}

\usepackage{epsfig}

\newtheorem{thm}{Theorem}
\newtheorem{pro}[thm]{Proposition}
\newtheorem{lem}[thm]{Lemma}
\newtheorem{cl}[thm]{Claim}
\newtheorem{cj}[thm]{Conjecture}

\newtheorem{df}[thm]{Definition}

\newcommand {\cbdo}{\hfill$\Box$}

\begin{document}
\begin{center}\textbf{{Further progress on Wojda's conjecture}}\\
\center{Maciej Cisi\'nski and Andrzej \.{Z}ak}\\
\end{center}
\begin{abstract}
Two digraphs of order $n$ are said to pack if they can be found as edge-disjoint subgraphs of the complete digraph of order $n$. It is well established that if the sum of the sizes of the two digraphs is at most $2n-2$, then they pack, with this bound being sharp. However, it is sufficient for the size of the smaller digraph to be only slightly below $n$ for the sum of their sizes to significantly exceed this threshold while still guaranteeing the existence of a packing.

In 1985, Wojda conjectured that for any $2 \leq m \leq n/2$, if one digraph has size at most $n - m$ and the other has size less than $2n - \lfloor n/m \rfloor$, then the two digraphs pack. It was previously known that this conjecture holds for $m = \Omega(\sqrt{n})$. In this paper,
we confirm it for $m \geq 93$ and $n \geq 31m$. 
\end{abstract}
\section{Introduction}
Given a digraph $D=(V, A)$, a set $W$ with $|W|\geq |V|$, and an injection $f: V \rightarrow W$, we define $f(D)$ as the digraph with vertex set $V(f(D))=W$ and arc set $A(f(D))$, where for any $u,v \in V$ an ordered pair $f(u)f(v) \in A(f(D))$ if and only if $uv \in A$.
\begin{df}
A packing of $D=(V, A)$ and $D' = (V', A')$ is a pair of injections 
$f: V \rightarrow W$ and $f': V' \rightarrow W$ such that $f(D)$ and $f'(D')$ are arc-disjoint.
\end{df}
When $|V| \leq |V'|=|W|$, which always holds in the sequel, we can simplify the notation by assuming $f'$ is the identity function and, without loss of generality, take $V'= W$. In this scenario, an injection 
$f: V \rightarrow V'$ is called a packing of $D=(V, A)$ and $D' = (V', A')$ if $f(D)$ and $D'$ are edge-disjoint. In other words, we say $D$ and $D'$ \emph{pack}. Furthermore, if $f(D)$ and $D'$ share at most $q$  arcs, $f$ is referred to as a $q$-\emph{near-packing}.

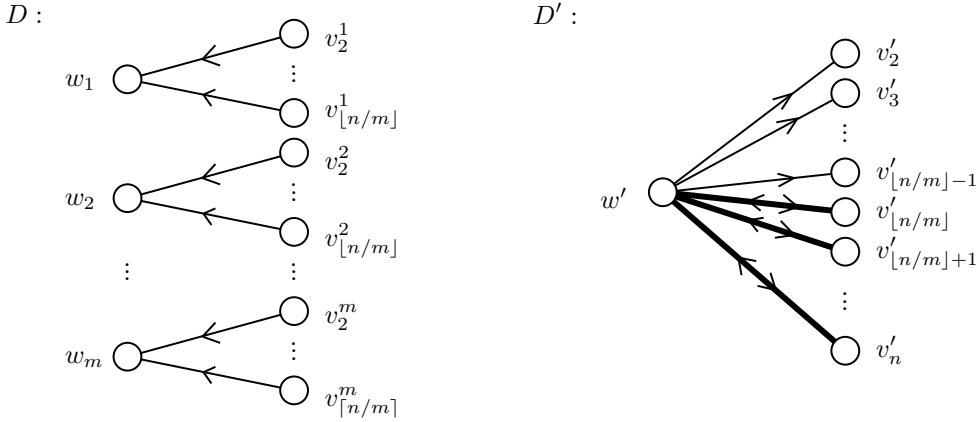
\begin{figure}[h]\label{fig}
    \centering

\tikzset{every picture/.style={line width=0.75pt}} 

\begin{tikzpicture}[x=0.75pt,y=0.75pt,yscale=-1,xscale=1]

\draw    (171.98,19) -- (87.98,42) ;
\draw    (172.01,59) -- (87.98,42) ;
\draw  [color={rgb, 255:red, 0; green, 0; blue, 0 }  ,draw opacity=1 ][fill={rgb, 255:red, 255; green, 255; blue, 255 }  ,fill opacity=1 ] (178.99,19.01) .. controls (179,15.14) and (175.86,12) .. (171.99,11.99) .. controls (168.12,11.99) and (164.98,15.12) .. (164.97,19) .. controls (164.97,22.87) and (168.1,26.01) .. (171.97,26.01) .. controls (175.85,26.02) and (178.99,22.88) .. (178.99,19.01) -- cycle ;
\draw  [color={rgb, 255:red, 0; green, 0; blue, 0 }  ,draw opacity=1 ][fill={rgb, 255:red, 255; green, 255; blue, 255 }  ,fill opacity=1 ] (94.99,42.01) .. controls (95,38.14) and (91.86,35) .. (87.99,34.99) .. controls (84.12,34.99) and (80.98,38.12) .. (80.97,42) .. controls (80.97,45.87) and (84.1,49.01) .. (87.97,49.01) .. controls (91.85,49.02) and (94.99,45.88) .. (94.99,42.01) -- cycle ;
\draw  [color={rgb, 255:red, 0; green, 0; blue, 0 }  ,draw opacity=1 ][fill={rgb, 255:red, 255; green, 255; blue, 255 }  ,fill opacity=1 ] (179.02,59.01) .. controls (179.02,55.14) and (175.89,52) .. (172.02,51.99) .. controls (168.15,51.99) and (165,55.12) .. (165,59) .. controls (165,62.87) and (168.13,66.01) .. (172,66.01) .. controls (175.87,66.02) and (179.02,62.88) .. (179.02,59.01) -- cycle ;
\draw    (357.98,99) -- (449.98,29) ;
\draw    (357.98,99) -- (449.99,49) ;
\draw    (357.98,99) -- (449.98,89) ;
\draw  [color={rgb, 255:red, 0; green, 0; blue, 0 }  ,draw opacity=1 ][fill={rgb, 255:red, 255; green, 255; blue, 255 }  ,fill opacity=1 ] (456.99,29.01) .. controls (457,25.14) and (453.86,22) .. (449.99,21.99) .. controls (446.12,21.99) and (442.98,25.12) .. (442.97,29) .. controls (442.97,32.87) and (446.1,36.01) .. (449.97,36.01) .. controls (453.85,36.02) and (456.99,32.88) .. (456.99,29.01) -- cycle ;
\draw  [color={rgb, 255:red, 0; green, 0; blue, 0 }  ,draw opacity=1 ][fill={rgb, 255:red, 255; green, 255; blue, 255 }  ,fill opacity=1 ] (456.99,49.01) .. controls (457,45.14) and (453.87,41.99) .. (450,41.99) .. controls (446.13,41.99) and (443,45.12) .. (442.99,48.99) .. controls (442.99,52.87) and (446.12,56.01) .. (449.99,56.01) .. controls (453.85,56.02) and (456.99,52.88) .. (456.99,49.01) -- cycle ;
\draw  [color={rgb, 255:red, 0; green, 0; blue, 0 }  ,draw opacity=1 ][fill={rgb, 255:red, 255; green, 255; blue, 255 }  ,fill opacity=1 ] (456.99,89.01) .. controls (457,85.14) and (453.86,82) .. (449.99,81.99) .. controls (446.12,81.99) and (442.98,85.12) .. (442.97,89) .. controls (442.97,92.87) and (446.1,96.01) .. (449.97,96.01) .. controls (453.85,96.02) and (456.99,92.88) .. (456.99,89.01) -- cycle ;
\draw [line width=2.25]    (357.98,99) -- (450.01,109.01) ;
\draw [line width=2.25]    (357.98,99) -- (450.01,129.01) ;
\draw [line width=2.25]    (357.98,99) -- (450.01,179) ;
\draw  [color={rgb, 255:red, 0; green, 0; blue, 0 }  ,draw opacity=1 ][fill={rgb, 255:red, 255; green, 255; blue, 255 }  ,fill opacity=1 ] (457.02,109.02) .. controls (457.02,105.15) and (453.89,102) .. (450.02,102) .. controls (446.15,102) and (443,105.13) .. (443,109) .. controls (443,112.87) and (446.13,116.02) .. (450,116.02) .. controls (453.87,116.02) and (457.02,112.89) .. (457.02,109.02) -- cycle ;
\draw  [color={rgb, 255:red, 0; green, 0; blue, 0 }  ,draw opacity=1 ][fill={rgb, 255:red, 255; green, 255; blue, 255 }  ,fill opacity=1 ] (457.02,129.02) .. controls (457.02,125.15) and (453.89,122) .. (450.02,122) .. controls (446.15,122) and (443,125.13) .. (443,129) .. controls (443,132.87) and (446.13,136.02) .. (450,136.02) .. controls (453.87,136.02) and (457.02,132.89) .. (457.02,129.02) -- cycle ;
\draw  [color={rgb, 255:red, 0; green, 0; blue, 0 }  ,draw opacity=1 ][fill={rgb, 255:red, 255; green, 255; blue, 255 }  ,fill opacity=1 ] (457.02,179.01) .. controls (457.02,175.14) and (453.89,172) .. (450.02,171.99) .. controls (446.15,171.99) and (443,175.12) .. (443,179) .. controls (443,182.87) and (446.13,186.01) .. (450,186.01) .. controls (453.87,186.02) and (457.02,182.88) .. (457.02,179.01) -- cycle ;
\draw  [color={rgb, 255:red, 0; green, 0; blue, 0 }  ,draw opacity=1 ][fill={rgb, 255:red, 255; green, 255; blue, 255 }  ,fill opacity=1 ] (364.99,99.01) .. controls (365,95.14) and (361.86,92) .. (357.99,91.99) .. controls (354.12,91.99) and (350.98,95.12) .. (350.97,99) .. controls (350.97,102.87) and (354.1,106.01) .. (357.97,106.01) .. controls (361.85,106.02) and (364.99,102.88) .. (364.99,99.01) -- cycle ;
\draw    (125.79,31.88) -- (134.73,35.07) ;
\draw    (125.79,31.88) -- (132.33,24.79) ;
\draw    (172,79) -- (88,102) ;
\draw    (125.81,91.88) -- (134.75,95.06) ;
\draw    (125.81,91.88) -- (132.35,84.78) ;
\draw    (172,159) -- (88,182) ;
\draw    (125.81,171.88) -- (134.75,175.06) ;
\draw    (125.81,171.88) -- (132.35,164.78) ;
\draw    (126.27,49.82) -- (130.33,55.52) ;
\draw    (126.27,49.82) -- (132.33,47.15) ;
\draw    (172.03,119) -- (88,102) ;
\draw    (126.29,109.82) -- (130.35,115.51) ;
\draw    (126.29,109.82) -- (132.35,107.15) ;
\draw    (172.03,199) -- (88,182) ;
\draw    (126.29,189.82) -- (130.35,195.51) ;
\draw    (126.29,189.82) -- (132.35,187.15) ;
\draw  [color={rgb, 255:red, 0; green, 0; blue, 0 }  ,draw opacity=1 ][fill={rgb, 255:red, 255; green, 255; blue, 255 }  ,fill opacity=1 ] (179.02,79.01) .. controls (179.02,75.14) and (175.89,72) .. (172.02,71.99) .. controls (168.15,71.99) and (165,75.12) .. (165,79) .. controls (165,82.87) and (168.13,86.01) .. (172,86.01) .. controls (175.87,86.02) and (179.02,82.88) .. (179.02,79.01) -- cycle ;
\draw  [color={rgb, 255:red, 0; green, 0; blue, 0 }  ,draw opacity=1 ][fill={rgb, 255:red, 255; green, 255; blue, 255 }  ,fill opacity=1 ] (179.02,119.02) .. controls (179.02,115.15) and (175.89,112) .. (172.02,112) .. controls (168.15,112) and (165,115.13) .. (165,119) .. controls (165,122.87) and (168.13,126.02) .. (172,126.02) .. controls (175.87,126.02) and (179.02,122.89) .. (179.02,119.02) -- cycle ;
\draw  [color={rgb, 255:red, 0; green, 0; blue, 0 }  ,draw opacity=1 ][fill={rgb, 255:red, 255; green, 255; blue, 255 }  ,fill opacity=1 ] (179.02,159.01) .. controls (179.02,155.14) and (175.89,152) .. (172.02,151.99) .. controls (168.15,151.99) and (165,155.12) .. (165,159) .. controls (165,162.87) and (168.13,166.01) .. (172,166.01) .. controls (175.87,166.02) and (179.02,162.88) .. (179.02,159.01) -- cycle ;
\draw  [color={rgb, 255:red, 0; green, 0; blue, 0 }  ,draw opacity=1 ][fill={rgb, 255:red, 255; green, 255; blue, 255 }  ,fill opacity=1 ] (179.02,199.01) .. controls (179.02,195.14) and (175.89,192) .. (172.02,191.99) .. controls (168.15,191.99) and (165,195.12) .. (165,199) .. controls (165,202.87) and (168.13,206.01) .. (172,206.01) .. controls (175.87,206.02) and (179.02,202.88) .. (179.02,199.01) -- cycle ;
\draw  [color={rgb, 255:red, 0; green, 0; blue, 0 }  ,draw opacity=1 ][fill={rgb, 255:red, 255; green, 255; blue, 255 }  ,fill opacity=1 ] (95.02,102.01) .. controls (95.02,98.14) and (91.89,95) .. (88.02,94.99) .. controls (84.15,94.99) and (81,98.12) .. (81,102) .. controls (81,105.87) and (84.13,109.01) .. (88,109.01) .. controls (91.87,109.02) and (95.02,105.88) .. (95.02,102.01) -- cycle ;
\draw  [color={rgb, 255:red, 0; green, 0; blue, 0 }  ,draw opacity=1 ][fill={rgb, 255:red, 255; green, 255; blue, 255 }  ,fill opacity=1 ] (94.99,182.01) .. controls (95,178.14) and (91.86,175) .. (87.99,174.99) .. controls (84.12,174.99) and (80.98,178.12) .. (80.97,182) .. controls (80.97,185.87) and (84.1,189.01) .. (87.97,189.01) .. controls (91.85,189.02) and (94.99,185.88) .. (94.99,182.01) -- cycle ;
\draw    (413.38,52.19) -- (423.1,49.33) -- (417.9,57.22) ;
\draw    (415.2,63.72) -- (425.9,61.74) -- (420.7,69.62) ;
\draw    (415.78,89.65) -- (424.15,91.82) -- (416.37,96.12) ;
\draw    (416.37,109.65) -- (426.04,106.41) -- (418.73,100.24) ;
\draw    (414.49,122) -- (424.27,120.53) -- (417.43,113.41) ;
\draw    (405.78,146.24) -- (414.98,148.41) -- (412.25,139.65) ;
\draw    (409.2,100.35) -- (401.21,103.47) -- (407.31,108.35) ;
\draw    (408.84,111.88) -- (401.21,113.12) -- (405.9,118.94) ;
\draw    (405.08,133.76) -- (396.04,131.94) -- (397.67,140.82) ;

\draw (55,38) node [anchor=north west][inner sep=0.75pt]   [align=left] {$\displaystyle w_{1}$};
\draw (55,98) node [anchor=north west][inner sep=0.75pt]   [align=left] {$\displaystyle w_2$};
\draw (55,178) node [anchor=north west][inner sep=0.75pt]   [align=left] {$\displaystyle w_{m}$};
\draw (174,132) node [anchor=north west][inner sep=0.75pt]  [rotate=-89.62] [align=left] {...};
\draw (186,12) node [anchor=north west][inner sep=0.75pt]   [align=left] {$\displaystyle v^1_{2}$};
\draw (186,48) node [anchor=north west][inner sep=0.75pt]   [align=left] {$\displaystyle v^1_{\lfloor n/m \rfloor}$};
\draw (186,74) node [anchor=north west][inner sep=0.75pt]   [align=left] {$\displaystyle v^2_2$};
\draw (186,113) node [anchor=north west][inner sep=0.75pt]   [align=left] {$\displaystyle v^2_{\lfloor n/m \rfloor}$};
\draw (186,154) node [anchor=north west][inner sep=0.75pt]   [align=left] {$\displaystyle v^m_2$};
\draw (186,196) node [anchor=north west][inner sep=0.75pt]   [align=left] {$\displaystyle v^m_{\lceil n/m \rceil}$};
\draw (25,3) node [anchor=north west][inner sep=0.75pt]   [align=left] {$\displaystyle D:$};
\draw (291,2) node [anchor=north west][inner sep=0.75pt]   [align=left] {$\displaystyle D':$};
\draw (452,147) node [anchor=north west][inner sep=0.75pt]  [rotate=-89.62] [align=left] {...};
\draw (464,20) node [anchor=north west][inner sep=0.75pt]   [align=left] {$\displaystyle v'_{2}$};
\draw (464,40) node [anchor=north west][inner sep=0.75pt]   [align=left] {$\displaystyle v'_{3}$};
\draw (464,100) node [anchor=north west][inner sep=0.75pt]   [align=left] {$\displaystyle v'_{\lfloor n/m\rfloor }$};
\draw (464,120) node [anchor=north west][inner sep=0.75pt]   [align=left] {$\displaystyle v'_{\lfloor n/m\rfloor +1}$};
\draw (464,170) node [anchor=north west][inner sep=0.75pt]   [align=left] {$\displaystyle v'_{n}$};
\draw (452,62) node [anchor=north west][inner sep=0.75pt]  [rotate=-89.62] [align=left] {...};
\draw (464,80) node [anchor=north west][inner sep=0.75pt]   [align=left] {$\displaystyle v'_{\lfloor n/m\rfloor -1}$};
\draw (325,95) node [anchor=north west][inner sep=0.75pt]   [align=left] {$\displaystyle w'$};
\draw (174,92) node [anchor=north west][inner sep=0.75pt]  [rotate=-89.62] [align=left] {...};
\draw (174,32) node [anchor=north west][inner sep=0.75pt]  [rotate=-89.62] [align=left] {...};
\draw (174,172) node [anchor=north west][inner sep=0.75pt]  [rotate=-89.62] [align=left] {...};
\draw (90,132) node [anchor=north west][inner sep=0.75pt]  [rotate=-89.62] [align=left] {...};

\end{tikzpicture}
\caption{Digraphs $D$ and $D'$ do not pack}
\end{figure}

In 1985, Wojda posed the following problem \cite{W}: 
for every $n, k, 1\leq k \leq n(n-1)$, determine the smallest number $\mu(n,k)$ such that there exist digraphs $D$ and $D'$ with $|A(D)|=k$ and $|A(D')|=\mu(n,k)$ for which there is no packing of $D$ and $D'$. Some known values of $\mu(n, k)$ are:  
\begin{align*}
\mu(n,1)=n(n-1),\; \mu(n,2)&={n \choose 2},\; \mu(n,3)={n \choose 2}-\lfloor \frac{n}{2} \rfloor  \text{ for } n\geq 7,  \\ 
\mu(n,n-1)&=n \text{ and } \mu(n,n)=n-1.
\end{align*}
For $0 < \alpha < 1/2$  the bound $\mu(n, \alpha n) = \Omega(n^{3/2})$ can be derived from a similar result concerning graphs \cite{Br}. Wojda proposed the following conjecture.
\begin{cj}\label{mainC}\cite{W}
For every $m$ satisfying $2\leq m \leq \frac{n}{2}$, 
\[\mu(n,n-m)=2n-\left\lfloor \frac{n}{m} \right\rfloor.\]
\end{cj}
The example provided in \cite{W} regarding two digraphs that cannot be packed indicates that $\mu(n,n-m)\leq 2n-\lfloor \frac{n}{m} \rfloor$. These digraphs, referred to as $D$ and $D'$, are formed through the following way. Define $n=a\lfloor n/m \rfloor + b\lceil n/m \rceil$ , where $a, b \in \mathbb{Z}$ with $a, b \geq 0$ and a+b=m. Thus, $V(D) = W\cup \bigcup_{i=1}^{a+b} V_i$ where $W = \{w_1,\dots, w_m\}$, and $V_i = \{v^i_2, \dots,v^i_{\lfloor n/m\rfloor} \}$ for $i = 1, \dots, a$, whereas $V_i = \{v^i_2, \dots,v^i_{\lceil n/m\rceil}\}$ for $i = a+1, \dots, b$. Additionally, $A(D) = \{vw_i: i=1, \dots, a+b, v \in V_i\}$. The digraph $D'$ has vertex set  $\{w', v'_2, \dots, v'_n\}$ and $2n-\lfloor \frac{n}{m} \rfloor$ arcs constructed as
\begin{align*}
&w'v'_j \text{ for } j = 2, \dots n, \\
&v'_jw' \text{ for } j = \lfloor n/m \rfloor, \dots n,
\end{align*}
see Figure \ref{fig}.

\newpage
These two digraphs do not pack, because the vertex $w'$ cannot be the image of any vertex from $V(D)$. Indeed, suppose first that $w'$ is an image of a vertex $w_i$ for some $i \in \{1,\dots, m\}$.  We then have to match the in-neighbors of $w_i$, i.e. the set $N^-_{D}(w_i)$, with the set $X'=V(D')\setminus N^-_{D'}[w']$. This, however, cannot be done because 
\[N^-_{D}(w_i) \geq \left\lfloor \frac{n}{m} \right\rfloor > \left\lfloor \frac{n}{m} \right\rfloor - 1 = |X'|.\] 
On the other hand, $w'$ cannot be the image of any $v^i_j$ either, as $d_{D'}^+(w') = n-1$ and $d_D^+(v_i^j)=1$.\\

The case $m=\frac{n}{2}$ of Conjecture \ref{mainC} follows directly from known results. The case $m = \Theta(n)$ can be proved relatively easily using another related result for graphs from \cite{BKN}. 
In \cite{KŻ}, Konarski and the second author proved Conjecture \ref{mainC} in cases where $m$ is relatively large.
\begin{thm}\label{KŻ}\cite{KŻ}
If $m\geq \sqrt{8n}+418275$, then
$$\mu(n,n-m)=2n-\left\lfloor \frac{n}{m}\right\rfloor.$$
\end{thm}

Our main result is the following theorem.

\begin{thm}\label{main}
For every $m\geq 93$ and $n\geq 31m$, \[\mu(n,n-m)=2n-\left\lfloor \frac{n}{m}\right\rfloor.\]
\end{thm}

We use conventional notation. For a digraph $D$, the vertex set is denoted by $V(D)$ and the arc set by $A(D)$. The indegree of a vertex $v$ in $D$ is represented by $d_{D}^{-}(v)$, while the outdegree is denoted by $d_{D}^{+}(v)$. The total degree, or simply the degree, of $v$ in $D$, noted as $d_{D}(v)$, is given by the formula $d_{D}(v)=d^{-}_{D}(v)+d^{+}_{D}(v)$. The neighborhood of a vertex $v$ is defined as $N^+_{D}(v)$ and $N^-_{D}(v)$, where
$$N^+_{D}(v)=\{u\in V(D): vu\in A(D)\} \text{ and } N^-_{D}(v)=\{u\in V(D): uv\in A(D)\}.$$
Additionally, we define $N_D(v)=N^+_{D}(v)\cup N^-_{D}(v)$. With this notation, one can also define both open and closed neighborhoods for sets of vertices. 
\begin{df} Let $X$ be a set of vertices of $D$, where $X=\{v_1,v_2,\dots,v_k\}$ for some vertices $v_1,v_2,\dots,v_k$ in $V(D)$. Then
$$N^+_D(X)=\bigcup_{v\in X}N^+_D(v)\setminus X\text{, } N^-_D(X)=\bigcup_{v\in X}N^-_D(v)\setminus X,$$  
$$N^+_D[X]=N^+_D(X)\cup X\text{ and } N^-_D[X]=N^-_D(X)\cup X.$$
\end{df}
\section{Preliminaries}
Rather than using induction, which would guarantee the existence of a packing for two digraphs by assuming that their subdigraphs pack, we will instead rely on the existence of a near-packing of the subdigraphs.

\begin{lem}\label{prod1}
Let $D, D'$ be digraphs on $n$ vertices. If $|A(D)|\cdot |A(D')|<(q+1)n(n-1)$, then there exists a $q$-near-packing of $D$ and $D'$.
\end{lem}

Proof. We shall use the probabilistic method. Consider the probability space whose $n!$ points are all the possible bijections $\sigma$ from $V(D)$ to $V(D')$, each with probability $\frac{1}{n!}$.
For any two arcs $i\in A(D)$ and $j\in A(D')$ we denote by $A_{ij}$ the \emph{undesirable} event that $j$ is an image of $i$. In other words,
$$A_{ij}=\{\sigma : V(D)\rightarrow V(D') \hspace{0.2cm}|\hspace{0.2cm} \sigma(i)=j\}.$$
Then $$P(A_{ij})=\frac{(n-2)!}{n!}=\frac{1}{n(n-1)}.$$
Let $X$ be a random variable with $$ X_{ij}=
\begin{cases}
1, \hspace{1cm}if \hspace{0.2cm} \sigma(i)=j,\\ 0, \hspace{1cm} \text{otherwise}.
\end{cases}$$
Then $$\mathbf{E}(\bigcup_{i\in A(D), j\in A(D')}X_{ij})= \sum_{i\in A(D), j\in A(D')}\mathbf{E}(X_{ij})= |A(D)|\cdot |A(D')|\cdot \frac{1}{n(n-1)}< q+1.$$ 
By applying the fundamental probabilistic technique, we can ascertain the existence of a permutation $\sigma^*$ that contains at most $q$ undesirable events $A_{ij}$. This implies there is a $q$-near-packing of $D$ and $D'$.\cbdo
\\

The next theorem, proved in \cite{BVW}, can also be derived from Lemma \ref{prod1} for $q=0$.
\begin{thm}\label{suma_arc}\cite{BVW}
 Let $D, D'$ be digraphs on $n$ vertices. If $|A(D)|+|A(D')| \leq 2n-2$, then $D$ and $D'$ pack.   
\end{thm}

The next lemma refines a well-known greedy algorithm for packing a tree with a graph or digraph, enhancing the likelihood of a successful packing by accounting for additional details.
\begin{lem}\label{lem1}
Let $F$ be an oriented forest that has at least $k+s$ components, $D'$ be a digraph such that $|V(F)| \leq |V(D')| = n$, and $S' \subset V(D')$ with $|S'| = k$. Let $V(D') = \{v'_1, \dots, v'_n \}$ with $d_{D'}(v'_i) \geq d_{D'}(v'_{i+1})$ for $i=1,\dots, n-1$. If

\begin{align*}
    a) \;\;\;&\max\left(d^-_{D'}(v'_i), d^+_{D'}(v'_i)\right)  \leq |V(D')| - i - s, \; i=1, \dots,s, \text{ and } \\
    &\left|N_{D'}(v'_{s+1}) \setminus \{v'_1, \dots, v'_s\}\right|  \leq |V(D')| - |V(F)|,
\end{align*}
or
\begin{align*}
    b) \;\;\;&\min\left(d^-_{D'}(v'_1), d^+_{D'}(v'_1)\right) \leq |V(D')| - |V(F)| \text{ and } \\
    &\left|N_{D'}(v'_2) \setminus \{v'_1\}\right| \leq |V(D')| - |V(F)|,
\end{align*}
then there is a packing $f$ of $F$ and $D'$ such that $\left(S' \cup \{v'_1, \dots, v'_s\}\right) \subset V(f(F))$.
\end{lem}
Proof. 
Let $S'' = \left|S' \cup \{v'_1, \dots, v'_s\}\right|$ with $|S''| = p \leq k+s$. Since $F$ has $k+s$ components and its underlying multigraph is a forest, there exists a set $L = \{l_1, \dots, l_{p}\}\subseteq V(F)$ such that $d_{F}(l_j)\leq 1$ and each $l_j$ belongs to a distinct component of $F$, for $j = 1, \dots, p$. We construct a packing $f$ in three stages. 

In stage 1 we select an appropriate number of vertices of $F$ and match them with all vertices of $D'$ that may have degree greater than $|V(D')| - |V(F)|$. Thus, in case a) we select $s$ vertices from $F$ and match them with the vertices of $\{v'_1, \dots, v'_s\}$ while in case b) we select only one vertex of $F$ and match it with $v'_1$. In both cases, we also match the neighbors of each selected vertex of $F$, if such neighbors exist.

Assume first that the conditions of case a) hold. Here, for each $i=1, \dots, s$, the neighbors of $f^{-1}(v'_i)$ will be matched with some vertices from $V(D')\setminus \{v'_1, \dots, v'_s\}$.  At the $i$-th step of stage 1, we set $f(l_i) = v'_i$. Since $l_i$ is the first matched vertex of its component of $F$, the packing property is preserved. 
If $d_{D'}(l_i) = 0$ we proceed to the next step $i+1$ of stage 1. Otherwise, let $u_i$ be the neighbor of $l_i$.  By the assumption on $\max(d^-_{D'}(v_i), d^+_{D'}(v_i))$, we have 
\begin{align*}
    &\min\left(|V(D')\setminus N^-_{D'}(v'_i)|, |V(D')\setminus N^+_{D'}(v'_i)|\right)  \geq i + s\text{ and so }\\
    &\min\left(|V(D')\setminus \left(N^-_{D'}(v'_i) \cup \{v'_1, \dots, v'_s\}\right)|, |V(D')\setminus \left(N^+_{D'}(v'_i) \cup \{v'_1, \dots, v'_s\}\right)|\right) \geq i.
\end{align*}
Since $l_i$ is the only matched neighbor of $u_i$, to complete step $i$, we need to find an unmatched vertex from $V(D')\setminus \left(N^-_{D'}(v'_i) \cup \{v'_1, \dots, v'_s\}\right)$ or $V(D') \setminus \left(N^+_{D'}(v'_i) \cup \{v'_1, \dots, v'_s\}\right)$ depending on whether $u_i$ is the in-neighbor or out-neighbor of $l_i$. At this moment, the number of matched vertices in $D'\setminus \{v'_1, \dots, v'_s\}$ is at most $i-1$. Therefore, there exists a vertex $u'$ that is still unmatched, allowing us to set $f(u_i) = u'$. 

Now, suppose that case b) holds. Without loss of generality, assume $d^+_{D'}(v'_1) \leq |V(D')|-|V(F)|$. Since $|A(F)| < |V(F)|$, there exists a vertex $u \in V(F)$ with $d^-_{F}(u) = 0$. We set $f(u)=v'_1$ and match each $v \in N^+_{F}(u)$ with an arbitrary vertex from $V(D')\setminus \left(N^+_{D'}(v'_1) \cup v'_1\}\right)$.

In Stage 2, we arbitrarily match the unmatched vertices of $L$ with the unmatched vertices of $S''$. If some vertices in $L$ remain unmatched after this step, we proceed by arbitrarily matching them with unmatched vertices of $D'$. During this stage, our primary focus is to ensure that each vertex in $S''$ is included in the image of $f$, as this is critical to satisfying the additional property of a packing.

In stage 3, we match the remaining unmatched vertices of $F$. At the start of this stage, in every component of $F$ there is an unmatched vertex that has exactly one matched neighbor or all vertices of this component are already matched. We select unmatched vertices from $F$ one by one, ensuring that the selected vertex has exactly one matched neighbor. This property is preserved throughout stage 3, as $F$ is acyclic and the matched vertices within any component of $F$ form a connected subgraph. 
Let $u$ be a yet unmatched vertex of $F$ and let $v$ be its matched neighbor. By the construction of stage 1, 
\begin{align*}
\left|N_{D'}(f(v)) \setminus \{v'_1, \dots, v'_s\}\right| \leq |V(D')| - |V(F)|.   
\end{align*} 
Furthermore, since the number of matched vertices in $D'$ is at most $|V(F)| - 1$ and $v'_1, \dots, v'_s$ are already matched, there exists at least one unmatched vertex $v' \in V(D')\setminus N_{D'}(f(v))$ allowing us to set $f(u) = v'$. We repeat this process for all unmatched vertices of $F$ until every vertex is matched. 

\cbdo

\begin{pro}\label{deg}
Let $D'$ be a digraph with $V(D') = \{v'_1, \dots, v'_n \}$ and $d_{D'}(v'_i) \geq d_{D'}(v'_{i+1})$ for $i=1,\dots, n-1$. Then for every $j\geq 2$
\begin{align*}
    d_{D'}(v'_j) &\leq \frac{|A(D')| + j(j-1)}{j}.
\end{align*}
\end{pro}
Proof. Note that 
\begin{align*}
jd_{D'}(v'_j)\leq d_{D'}(v'_1) + \dots + d_{D'}(v'_j)\leq |A(D')| + j(j-1),    
\end{align*}
which proves the statement of Proposition \ref{deg}.
\cbdo 
\section{Proof of Theorem \ref{main}}
We assume that $V(D') = \{v'_1, \dots, v'_n \}$ with $d_{D'}(v'_i) \geq d_{D'}(v'_{i+1})$ for $i=1,\dots, n-1$, $|A(D')| = 2n - \lfloor n/m \rfloor -1$ and $|A(D)| = n-m$. We will show that $D$ and $D'$ pack. 
Let $\gamma'$ be the number of vertices of degree at most $3$ in $D'$. Then
\begin{align*}
   4(n-\gamma') \leq \sum_{v \in V(D')} d_{D'}(v) < 4n - \frac{2n}{m}. 
\end{align*}
Hence, 
\begin{align}\label{n3}
    \gamma' > \frac{n}{2m} > 15,
\end{align}
by the assumption on $n$. Since $|A(D)| = n - m$, $D = T_1 \cup \dots \cup T_m \cup R$, where $T_i$ are oriented trees. Without loss of generality, we assume that $|T_i| \leq |T_{i+1}|$, $i = 1, \dots, m-1$.
Let $F_k = T_1 \cup \dots \cup T_k$. 
\begin{cl}\label{cl1}
For every $k$, $|F_k| \leq k\frac{n - |R|}{m}$.
\end{cl}
Proof of Claim \ref{cl1}. Since $|T_j| \geq |T_i|$ for $j > i$, we have
\begin{align*}
  \sum_{j=1}^k \frac{m}{k}|T_j| &\leq \sum_{j=1}^{k-1} |T_j| + (m-k+1)|T_k|
  \leq \sum_{j=1}^{m} |T_j| = n -|R|,
\end{align*}
which proves the claim. \cbdo
\\

By \eqref{n3}, we can choose four independent vertices, say $u'_1, \dots, u'_4$ of $D'$, each of degree at most 3. Let $D'' = D'-\{u'_1, \dots, u'_4\}$, $S' = N_{D'}(\{u'_1, \dots, u'_4\})$ and $T' = N_{D'}[\{u'_1, \dots, u'_4\}]$. Clearly $V(D'') = n - 4$, $|S'| \leq 12$ and $|T'| \leq 16$.  
We divide our proof into two cases according to the value of $d_{D'}(v'_1)$. 

Suppose first that $d_{D'}(v'_1) \leq n - 14$. We will use Lemma \ref{lem1} with $s = 5$ ($s=2$ would be good as well, but $s=5$ minimizes the value of $m$) and $k=12$ for $D''$ and $F_{17}$. Let $V(D'') = \{v''_1, \dots, v''_{n-4} \}$ with $d_{D''}(v''_i) \geq d_{D''}(v''_{i+1})$ for $i=1,\dots, |V(D'')|-1$. Clearly 
\begin{align*}
\max\left( d^-_{D''}(v''_1), d^+_{D''}(v''_1), d^-_{D''}(v''_2), d^+_{D''}(v''_2) \right) \leq d_{D'}(v'_1) \leq n - 14 \leq |V(D'')| - 10.   
\end{align*} 
Furthermore, by Proposition \ref{deg} and Claim \ref{cl1}

\begin{align*}
    d_{D''}(v''_6) \leq \frac{2n-\lfloor n/m \rfloor - 1 + 30}{6} \leq n-4 - \frac{17n}{m},
\end{align*}
because $m \geq 26$. Thus, by Lemma \ref{lem1}, there exists a packing $\phi$ of $D''$ and $F_{17}$ such that $\left( S' \cup \{v'_1, \dots,  v'_5 \} \right) \subset \phi(V(F_{17}))$. Let $l_{18} \in T_{18}$ and $l_{19} \in T_{19}$. Let $u_{18}$, $u_{19}$ be the neighbors of $l_{18}$ and $l_{19}$, respectively. 
Let $H' = D' - \left(\phi(V(F_{17}))\cup \{u'_1, \dots, u'_4\}\right)$ and $H = D - \left(V(F_{17})\cup \{l_{18}, l_{19}, u_{18}, u_{19}\} \right)$. 
By Claim \ref{cl1} and since $m$ is relatively large
\begin{align*}
    |V(H)| = |V(H')| &\geq n - \frac{17n}{m} - 4 \\ 
    |A(H)|\cdot |A(H')| &< n\left(2n-\frac{n}{m}\right) \leq 3\left(n - \frac{17n}{m} - 4\right)\left(n - \frac{17n}{m} - 5\right).
\end{align*}
Thus, by Lemma \ref{prod1} there is a 2-near packing $\psi$ of $H$ and $H'$. Let $x', y'$ contain the vertex cover of $\psi(H) \cap H'$. Then 
\begin{align*}
    f(v) =
\begin{cases}
    u'_1 & \text{ if } v = \psi^{-1}(x'), \\
    u'_2 & \text{ if } v = \psi^{-1}(y'), \\
    u'_3 & \text{ if } v = u_{18},\\
    u'_4 & \text{ if } v = u_{19},\\
    x' & \text{ if } v = l_{18},\\
    y' & \text{ if } v = l_{19},\\
    (\phi \cup \psi)(v) & \text{ otherwise}, 
\end{cases}
\end{align*}
is a packing of $D$ and $D'$.

So, we may assume that $d_{D'}(v'_1) \geq n - 13$. 
We will examine four distinct subcases. For three of these, Lemma \ref{lem1} will be applied to demonstrate that there is a packing $\phi$ of $T_1$ and $D'$, such that $v'_1 \in \phi(V(T_1))$. Subsequently, $\phi$ will be expanded for all three scenarios to create a packing of $D$ and $D'$. In the final subcase, we will develop a packing $\psi$ of $F_2$ and $D'$ such that $v'_1, v'_2 \in \psi(V(F_2))$, and proceed to extend $\psi$ into a packing of $D$ and $D'$. Without loss of generality, assume $d^-_{D'}(v'_1) = \max(d^-_{D'}(v'_1), d^+_{D'}(v'_1))$. 

\textbf{Subcase a)} If $d^-_{D'}(v'_1) = n-1$ then 
\begin{align*}
d^+_{D'}(v'_1) \leq 2n-\lfloor n/m \rfloor -1 - (n-1) = n - \lfloor n/m \rfloor \leq |V(D')|-|V(T_1)|,
\end{align*}
by Claim \ref{cl1}. Furthermore,
\begin{align}\label{nv2}
     |N_{D'}(v'_2)\setminus\{v'_1\}| &\leq |A(D')| - d_{D'}(v'_1) \leq 2n - \lfloor n/m \rfloor -1 - (n-1) \\
     &\leq n - \lfloor n/m \rfloor \leq |V(D')|-|V(T_1)|. \nonumber
\end{align}
Thus, conditions b) of Lemma \ref{lem1} hold and so there exists a packing $\phi$ of $D'$ and $T_1$ such that $v'_1 \in \phi(V(T_1)$.

\textbf{Subcase b)} If $d^-_{D'}(v'_1) = n-2$ then we can assume that $d_{D'}(v'_1) \geq n-1$. Indeed, if $d^+_{D'}(v'_1) = 0$ let $H' = D' - v'_1$ and $H = D-u$ where $d^-_{D}(u) = 0$ (such a vertex $u$ exists because $|A(D)| < n$).  Then 
\begin{align*}
 |A(H)|+|A(H')| \leq n-m + 2n - \lfloor n/m \rfloor -1 - (n-2) \leq 2n - \lfloor n/m \rfloor + 1 \leq 2(n-1) - 2, 
\end{align*}
by the assumption on $n$. Hence by Theorem \ref{suma_arc} there is a packing $h$ of $H$ and $H'$. Since $d^+_{D'}(v'_1) = 0$ and $d^-_{D}(u) = 0$, $f$ such that $f(u) = v'_1$ and $f(v) = h(v)$ for each $v \neq u$ is a packing of $D$ and $D'$. Thus $d^+_{D'}(v'_1) \geq 1$ and so $d_{D'}(v'_1) \geq n-1$.
Therefore, \eqref{nv2} holds and so conditions a) of Lemma \ref{lem1} are satisfied for $s=1$ and $k=0$. Thus, there is a packing as required of $T_1$ and $D'$. 

\textbf{Subcase c)} If $d^-_{D'}(v'_1) \leq n-3$ and $|N_{D'}(v'_2)\setminus\{v'_1\}|\leq n - \lfloor n/m \rfloor$ then conditions a) of Lemma \ref{lem1} hold and so there exists a packing as required of $T_1$ and $D'$.

In all three subcases let $\phi$ be a packing of $T_1$ and $D'$ such that $v'_1 \in \phi(V(T_1)$. Let $H' = D' - \phi(V(T_1))$ and $H = D - V(T_1)$. Thus, $|A(H)|\leq n - m - |V(T_1)|+1$, $|V(H)|= |V(H')| = n- |V(T_1)|$ and $|A(H')|\leq 2n - \lfloor n/m \rfloor - 1 - (n - 13) = n-\lfloor n/m \rfloor + 12$. Then
\begin{align*}
    |A(H)|+|A(H')| & \leq 2n - m - \lfloor n/m \rfloor - |V(T_1)|+13 \\
    & = 2n-2 - \left(\lfloor n/m \rfloor + m - 15 + |V(T_1)| \right) \\
    &\leq 2\left(n - |V(T_1)|\right) - 2 =  2|V(H)| - 2, 
\end{align*}
since $|V(T_1)| \leq \lfloor n/m \rfloor$ and $m \geq 15$. Hence, $H$ and $H'$ pack by Theorem \ref{suma_arc}, and so $D$ and $D'$ pack as well. 

\textbf{Subcase d)} Finally, assume that $d_{D'}(v'_1) \geq n - 13$, $d^-_{D'}(v'_1) \leq  n - 3$ and $|N_{D'}(v'_2)\setminus\{v'_1\}| > n - \lfloor n/m \rfloor$. 
Then 
\begin{align*}
   |N_{D'}(v'_2)\setminus\{v'_1\}| &\leq |A(D')| - d_{D'}(v'_1) \leq 2n - \lfloor n/m \rfloor -1 - (n-13) \\
     &\leq n - \lfloor n/m \rfloor + 12\leq |V(D')|-4, \nonumber 
\end{align*}
by the assumption on $n$. 
Furthermore, by Proposition \ref{deg} and Claim \ref{cl1}, 
\begin{align*}
    d_{D'}(v'_3) \leq \frac{2n-n/m - 1 + 6}{3} \leq n - \frac{2n}{m} \leq n - |V(F_2)|.
\end{align*}
Hence, by Lemma \ref{lem1} with $s=2$ and $k=0$, there exists a packing $\psi$ of $D'$ and $F_2$ such that $\{v'_1, v'_2\} \subset \psi(V(F_2))$. 
Let $H' = D' - \psi(V(F_2))$ and $H = D - V(F_2)$. Thus, $|A(H)|\leq n - m - |V(F_2)|+1$, $|V(H)|= |V(H')| = n - |V(F_2)|$ and 
\begin{align*}
|A(H')|\leq 2n - \lfloor n/m \rfloor -1- (n - 13) - (n - \lfloor n/m \rfloor) = 12. 
\end{align*}
Thus,
\begin{align*}
  |A(H)|+|A(H')|  \leq n - m - |V(F_2)|+ 13 \leq (n -|V(F_2)|) - 2 \leq 2(n -|V(F_2)|) - 2, 
\end{align*}
since $m \geq 15$. Hence, $H$ and $H'$ pack by Theorem \ref{suma_arc}, and so $D$ and $D'$ pack as well.  
\cbdo

\end{document}